\documentclass[11pt,a4paper]{article}

\usepackage[utf8]{inputenc}
\usepackage[english]{babel}
\usepackage{amsmath,amssymb,amsfonts,amsthm}
\usepackage{geometry}
\numberwithin{equation}{section}

\newtheorem{theorem}{Theorem}[section]
\newtheorem{lemma}[theorem]{Lemma}
\theoremstyle{definition}
\newtheorem{remark}[theorem]{Remark}

\newcommand{\R}{\mathbb{R}}
\DeclareMathOperator{\divop}{div}
\DeclareMathOperator{\limsupop}{lim\,sup}

\title{\textbf{Bernstein's theorem for variational integrals\\ of linear growth and radial structure}}
\author{\textbf{Martin Fuchs}\quad \& \quad \textbf{Michael Bildhauer}}
\date{}

\begin{document}

\maketitle

\begin{abstract}
We consider entire solutions $u: \R^2 \rightarrow \R$ of the Euler-Lagrange equation associated to the variational integral $\int_{\Omega} g(|\nabla u|)\,dx$ with a strictly convex density $g: [0,\infty)\rightarrow \R$ being of linear growth. We show that the condition $\int_{0}^{\infty} t\,g''(t)\,dt < \infty$ implies the Bernstein property, which means that $u$ must be an affine function. If this condition on g is weakened, we still have some partial Bernstein results.
\vspace{0.3cm}

\end{abstract}

\section{Introduction}

The famous theorem of Bernstein for nonparametric minimal surfaces states that an entire solution $u \in C^2(\R^2)$ of the minimal surface equation
\begin{equation}
\divop \left( \frac{\nabla u}{\sqrt{1+|\nabla u|^2}} \right) = 0
\end{equation}
must be an affine function, which means in more geometrical terms: if the mean curvature of the surface $S := \{(x_1, x_2, u(x_1, x_2)) : (x_1, x_2) \in \R^2\}$ vanishes, then $S$ is just a plane. The proof can be traced in Bernstein's work [1] dated in the period from 1915 to 1917, and a complete proof together with some extensions to equations of the type (1.1) is presented in [2].

Observing that (1.1) represents the Euler-Lagrange equation for local minimizers $u: \R^2 \rightarrow \R$ of the area integral $\int_{\Omega} g_0(|\nabla u|)\,dx$, $g_0(t) := \sqrt{1+t^2}$, defined on domains $\Omega \subset \R^2$, we ask the question of reasonable conditions to be imposed on a density function $g: [0,\infty) \rightarrow \R$ under which the Bernstein property holds for the functional $\int_{\Omega} g(|\nabla u|)\,dx$ being of radial structure.\\ More precisely, this means: are the affine functions $u: \R^2 \rightarrow \R$ the only entire solutions of the equation
\begin{equation}
\divop \left( \frac{g'(|\nabla u|)}{|\nabla u|} \nabla u \right) = 0
\end{equation}
under suitable assumptions on g? A positive answer to the validity of the Bernstein property has been given in the different setting of parametric integrands for instance by Bers [3], Finn [4], Jenkins [5] and Simon [6], who discuss the Euler - Lagrange equations for energies of the form $\int_{\Omega} G((1,\nabla u))\,dx$, $\Omega \subset \R^2$, with density $G: \R^3 \rightarrow \R$ of linear growth being positively homogeneous of degree 1 and complemented with appropriate ellipticity conditions. In conclusion, for densities $g$ taking the form $g(t) = \tilde{g}(\sqrt{1+t^2})$ with $\tilde{g}$ of linear growth and a sufficient degree of regularity, the Bernstein problem is solved.

Let us now give a precise formulation of our assumptions concerning the function g occurring in (1.2). Throughout this paper we consider functions $g: [0,\infty) \rightarrow \R$ at least of class $C^2$ such that
\begin{equation}
g'(0) = 0, \quad g''(t) > 0 \text{ for } t \geqslant 0,
\end{equation}
\begin{equation}
m(t-1) \leqslant g(t) \leqslant M(t+1) \text{ on } [0,\infty),
\end{equation}
where $m, M$ denote positive constants. Moreover, we require for the moment the upper bound
\begin{equation}
g''(t) \leqslant c(1+t)^{-\mu}, \quad t \geqslant 0,
\end{equation}
with exponent $\mu > 1$ and for a constant $c > 0$.

Note that the minimal surface density $g_0(t) = \sqrt{1+t^2}$ fulfills (1.3)--(1.5) with $\mu = 3$ as the optimal choice for the exponent $\mu$ occurring in the upper bound for $g''$. It turns out that the size of the exponent $\mu$ from (1.5) plays a key role in the decision of the Bernstein property for equation (1.2). To our knowledge, the first important contribution is contained in Theorem 1.4 from the work [7] of Farina, Sciunzi and Valdinoci, which we rephrased in Theorem 1.2 of [8] with the result that $\mu \geqslant 3$ in (1.5) is a sufficient condition for the Bernstein property of equation (1.2). In order to weaken the restriction $\mu \geqslant 3$ we discuss the assumption
\begin{equation}
\int_{0}^{\infty} s\,g''(s)\,ds < \infty,
\end{equation}
which follows from (1.5) in the case that
\begin{equation}
\mu > 2.
\end{equation}

Condition (1.6) has a geometric as well as an analytic meaning: in [9] we introduced the notion of $g$-catenoids for densities $g$ with (1.3)--(1.5). It turns out (see Lemma 2 in [9]) that these radially symmetric solutions defined for $|x| > 1$ show the same behaviour as the standard catenoid corresponding to the function $g_0(t) = \sqrt{1+t^2}$, whereas the failure of (1.6) generates a completely different behaviour of these functions. So we venture the guess that this reversal of the geometric properties of radially symmetric solutions of (1.2) suggests the breakdown of the Bernstein property, if (1.6) is violated. At the same time condition (1.6) is of essential importance for the investigations carried out in [10]: using (1.6) the concept of nonparametric $\mu$-surfaces is introduced, whose analytic description is very similar to the minimal surface case. For details we refer to Theorem 1.3 of [10]. Another indication showing the importance of (1.6) is the existence result of Beck, Bulíček and Maringová [11], who showed that the violation of (1.6) guarantees the existence of a unique solution to the Dirichlet boundary value problem for equation (1.2) on rather general -- not necessarily convex -- bounded domains $\Omega \subset \R^n, n \geqslant 2$. In other words: if (1.6) does not hold, then we have the same results as in the standard elliptic setting, which again suggests the failure of the Bernstein property.

In order to prove the Bernstein theorem under the assumption (1.6) we benefit essentially from the work [12] of Mikljukov, who investigated the behaviour of the solutions to certain types of equations on surfaces $S$ in $\R^3$ in terms of the geometric structure of $S$. From his Theorem 5.1 we obtain by just checking the assumptions from Section 5.2 in [12]

\begin{theorem} \label{thm1.1}
Suppose that the function $g: [0,\infty) \rightarrow \R$ satisfies the conditions (1.3)--(1.6) and assume in addition that
\begin{equation}
g \in C^3([0,\infty)), \quad g^{(3)}(0) = 0.
\end{equation}
Then we have the Bernstein property for equation (1.2): if $u \in C^3(\R^2)$ is an entire solution of (1.2), then $u$ is an affine function.
\end{theorem}

\begin{remark}
The assumption that $u$ and $g$ are of class $C^3$ enters through the methods used in the proofs of Theorems III and 4.2 of [12].
\end{remark}

\begin{remark} \label{rem1.3}
We strongly believe that the Bernstein property fails, if (1.6) or (1.7) are violated. However, a counterexample is still missing.
\end{remark}

Regarding the last remark it seems to be of interest to single out a class of entire solutions of equation (1.2) for which the Bernstein property holds just under the assumptions (1.3)--(1.5) imposed on the density function $g$. For this purpose we introduce the radial and tangential derivative of a function $u: \R^2 \rightarrow \R$ through the formulas
\[
\partial_r u(x) := \nabla u(x) \cdot \frac{x}{|x|}, \quad \partial_t u(x) := \nabla u(x) \cdot \frac{x^\perp}{|x|}
\]
for points $x \in \R^2 - \{0\}$, where we have set $x^\perp := (-x_2, x_1)$. It clearly holds
\[
\nabla u(x) = \partial_r u(x) \frac{x}{|x|} + \partial_t u(x) \frac{x^\perp}{|x|},
\]
and we have the following result:

\begin{theorem}
Suppose that $g \in C^2([0,\infty))$ satisfies the conditions (1.3) and (1.4), where in (1.3) the case $g''(0) = 0$ is admissible. We further assume that
\begin{equation}
\limsupop_{t \rightarrow \infty} \left[ g''(t) \, t \, \ln^2(t) \right] < \infty.
\end{equation}
If $u \in C^2(\R^2)$ is an entire solution of equation (1.2) with the property
\begin{equation}
\ln(1+|\nabla u|)\,|\partial_t u| \leqslant K(1+|\nabla u|) \quad \text{on } \R^2
\end{equation}
for a finite constant $K$, then $u$ must be affine.
\end{theorem}

\begin{remark} \label{rem1.6}
Condition (1.10) clearly holds in the case that $\partial_t u \in L^\infty(\R^2)$ or if we assume that $|\partial_t u| \leqslant c(|\partial_r u|^{1-\varepsilon} + 1)$ for some arbitrarily small number $\varepsilon > 0$ and a constant $c \in (0,\infty)$. More generally, we can replace (1.10) by the inequality $|\partial_t u| \leqslant c(1+|\partial_r u|)/\ln(e+|\partial_r u|)$.
\end{remark}

\begin{remark}
Inequality (1.9) is a weak form of assumption (1.5) in which $\mu > 1$ is required.
\end{remark}
\noindent Our paper is organized as follows: in Section 2 we give the proof of Theorem \ref{thm1.1} based on Mikljukov's work [12], the proof of Theorem 1.4 is presented in Section 3.

\section{Proof of Theorem 1.1}

Let the assumptions of Theorem \ref{thm1.1} hold and consider an entire solution $u \in C^3(\R^2)$ of equation (1.2). For $p \in \R^2$ we let $f(p) := g(|p|)$ and $a(p) := \nabla f(p) = \frac{g'(|p|)}{|p|}p$, where $a(0) = 0$ on account of $g'(0) = 0$. Note that $f \in C^2(\R^2)$ and in consequence $a \in C^1(\R^2, \R^2)$ holds in view of our assumption $g \in C^2([0,\infty))$ together with $g'(0) = 0$. Equation (1.2) then takes the form
\begin{equation}
\sum_{i=1}^{2} \frac{\partial}{\partial x_i} a_i (\nabla u) = 0 ,
\end{equation}
and in accordance with [12], Section 5.2, we have to check ($a_{ij} := \frac{\partial a_i}{\partial p_j} = \frac{\partial^2 f}{\partial p_i \partial p_j}$, $i,j = 1,2$):
\begin{enumerate}
\item[i)] For a constant $k \in (0,\infty)$ and with a suitable function $\Theta: [0,\infty) \rightarrow [0,\infty)$ such that
\begin{equation}
\Theta \in C^1([0,\infty)), \quad \Theta'(t) \leqslant 0, \quad \int_{0}^{\infty} \Theta(t)dt < \infty
\end{equation}
it holds
\begin{equation}
(1+p_1^2)a_{11}(p) + p_1 p_2(a_{12}(p) + a_{21}(p)) + (1+p_2^2)a_{22}(p) \leqslant k  \Theta(|p|) \sqrt{1+|p|^2}
\end{equation}
for any $p = (p_1, p_2) \in \R^2$.
\item[ii)] With a finite constant $c>1$ we have 
\begin{equation}
4 a_{11}(p) a_{22}(p) - (a_{12}(p) + a_{21}(p))^2>\frac{c}{4} (a_{12}(p) - a_{21}(p))^2
\end{equation}
for all points $p \in \R^2$.
\item[iii)] The components $a_i$, $i=1,2$, of the field $a: \R^2 \rightarrow \R^2$ are twice continuously differentiable functions.
\end{enumerate}
Accepting this for a moment and noting that conditions (2.1) and (2.2) from [12] trivially hold, Theorem 5.1 in [12] shows that $\nabla u(x)$ has a limit $\overline{p} \in \R^2$ as $x \rightarrow \infty$, hence $\nabla u$ is a globally bounded function, and from Moser's work [13] we deduce our claim.\\\\
\noindent \textbf{ad i):} We observe the identity ($i,j=1,2$)
\begin{equation}
a_{ij}(p) = \frac{g'(|p|)}{|p|} \delta_{ij} + \frac{p_i p_j}{|p|^2} \left[ g''(|p|) - \frac{g'(|p|)}{|p|} \right]
\end{equation}
valid for $p \in \R^2$ and get for the left-hand side of (2.3) by using (2.5)
\begin{align*}
(1+p_1^2)a_{11}(p) &+ p_1 p_2 \bigl( a_{12}(p) + a_{21}(p) \bigr) + (1+p_2^2)a_{22}(p) \\
&= a_{11}(p) + a_{22}(p) + \sum_{i,j=1}^{2} a_{ij}(p) p_i p_j \\
&= 2\frac{g'(|p|)}{|p|} + \frac{p_1^2 + p_2^2}{|p|^2} \left[ g''(|p|) - \frac{g'(|p|)}{|p|} \right] \\
&\quad + \sum_{i,j=1}^{2} a_{ij}(p) p_i p_j.
\end{align*}
Recalling $g'(0) = 0$ and noting that (1.4) gives the boundedness of $g'$, we find ($t := |p|$)\\\\
\(g^{\prime }(t)/t\leqslant c_{1}(1+t^{2})^{-1}(1+t^{2})^{1/2},\)\\\\
and for $i = 1,2$ it follows from (1.5)\\\\
\(\frac{p_{i}p_{i}}{t^{2}}g^{\prime \prime }(t)\leqslant c_{2}(1+t^{2})^{1/2}(1+t^{2})^{-\frac{\mu }{2}-\frac{1}{2}},\)\\\\
where we use the symbol $c_l$ to denote constants $\in (1,\infty)$ independent of $p$. We emphasize that in (1.5) we just require that $\mu > 1$, hence we obtain:
\begin{equation}
\text{l.h.s. of (2.3)} \leqslant c_3 (1+t^2)^{-1} \sqrt{1+t^2} + \sum_{i,j=1}^{2} a_{ij}(p) p_i p_j.
\end{equation}
It further holds (again by (2.5))
\\\(\begin{aligned}\sum _{i,j=1}^{2}a_{ij}(p)p_{i}p_{j}&=g^{\prime }(|p|)\frac{1}{|p|}p_{1}^{2}+\frac{p_{1}^{4}}{|p|^{2}}[g^{\prime \prime }(|p|)-g^{\prime }(|p|)/|p|]+g^{\prime }(|p|)\frac{1}{|p|}p_{2}^{2}\\ &\quad +\frac{p_{2}^{4}}{|p|^{2}}[g^{\prime \prime }(|p|)-g^{\prime }(|p|)/|p|]+2\frac{p_{1}^{2}p_{2}^{2}}{|p|^{2}}[g^{\prime \prime }(|p|)-g^{\prime }(|p|)/|p|]=\\ &=g^{\prime }(|p|)|p|^{-3}[p_{1}^{2}|p|^{2}-p_{1}^{4}+p_{2}^{2}|p|^{2}-p_{2}^{4}-2p_{1}^{2}p_{2}^{2}]\\ &\quad +g^{\prime \prime }(|p|)|p|^{-2}[p_{1}^{4}+p_{2}^{4}+2p_{1}^{2}p_{2}^{2}]=\\ &=t^{2}g^{\prime \prime }(t),\quad t:=|p|.\end{aligned}\)\\\
From (2.6) we get for the left-hand side of (2.3)
\[
\text{lhs of (2.3)} \le c_4 ((1 + t^2)^{-1}+t g''(t)) \sqrt{1 + t^2}
\]
and if we let $\Theta(t) = (1 + t^2)^{-1}+t g''(t)$, then we obtain 
(2.3) on account of our assumption (1.6).\\\\
\noindent \textbf{ad ii):} Due to the structure of the coefficients $a_{ij}(p)$ inequality (2.4) is equivalent to the validity of
\begin{equation}
a_{11}(p) a_{22}(p) > a_{12}(p)^2 \quad \text{for all } p \in \R^2.
\end{equation}
It holds\\

\begin{equation}
(a_{12}(p)^{2}=|p|^{-4}p_{1}^{2}p_{2}^{2}[g^{\prime \prime }(|p|)-g^{\prime }(|p|)/|p|]^{2})
\end{equation}
and

\(\begin{aligned}\text{l.h.s.\ of\ (2.7)}&=\left(\frac{g^{\prime }(|p|)}{|p|}+\frac{p_{1}^{2}}{|p|^{2}}\left[g^{\prime \prime }(|p|)-\frac{g^{\prime }(|p|)}{|p|}\right]\right)\left(\frac{g^{\prime }(|p|)}{|p|}+\frac{p_{2}^{2}}{|p|^{2}}\left[g^{\prime \prime }(|p|)-\frac{g^{\prime }(|p|)}{|p|}\right]\right)\\ &\quad =|p|^{-4}p_{1}^{2}p_{2}^{2}\left[g^{\prime \prime }(|p|)-\frac{g^{\prime }(|p|)}{|p|}\right]^{2}\\ &+\left(\frac{g^{\prime }(|p|)}{|p|}\right)^{2}+\frac{g^{\prime }(|p|)}{|p|}\left[g^{\prime \prime }(|p|)-\frac{g^{\prime }(|p|)}{|p|}\right]\\ &=\frac{1}{|p|}g^{\prime }(|p|)g^{\prime \prime }(|p|)+|p|^{-4}p_{1}^{2}p_{2}^{2}\left[g^{\prime \prime }(|p|)-\frac{g^{\prime }(|p|)}{|p|}\right]^{2}\\.\end{aligned}\)\\\\
Thus we obtain (2.7) as a consequence of assumption (1.3).\\\\
\noindent \textbf{ad iii):} We return to equation (2.5) to get for $k = 1,2$ and $p \in \R^2$, $p \neq 0$ ($t := |p|$)

\begin{eqnarray*}
\partial_k a_{ij}(p) &= &\frac{d}{dt}\left(\frac{g'(t)}{t}\right) \delta_{ij} \frac{p_k}{|p|} + \partial_k \left(\frac{p_i p_j}{|p|^2}\right)[g''(t)-g'(t)/t] \\
&+& \frac{p_i p_j}{|p|^2}\frac{d}{dt}[g''(t)-g'(t)/t]\frac{p_k}{|p|} =: (\alpha) + (\beta) + (\gamma).
\end{eqnarray*}

We calculate the single terms for the limit $t \rightarrow 0$: using a mid-value $\tilde{t} \in (0,t)$ based on $g'(0) = 0$, it holds (passing to the limit $t\to 0$):
\begin{align*}
\lvert (\alpha) \rvert &\leqslant \frac{\lvert t g''(t) - g'(t) \rvert}{t^2} = \frac{\lvert g''(t) - g''(\tilde{t}) \rvert}{t} \leqslant \max_{[0,t]} \lvert g^{(3)} \rvert \rightarrow 0), \\
\lvert (\beta) \rvert &\leqslant c \frac{1}{t} \lvert g''(t) - g'(t)/t \rvert \rightarrow 0,
\\\
\lvert (\gamma)\rvert&\leqslant\left\rvert \frac{d}{dt}(g''(t) - g'(t)/t) \right\rvert = \left\vert g^{(3)}(t)-\frac{t g''(t) - g'(t)}{t^2} \right\lvert = \left\rvert g^{(3)}(t) - \frac{1}{t}(g''(t) - g''(\tilde{t})) \right\rvert \rightarrow 0 .
\end{align*}
To sum up, it is proved that $\lim_{p \rightarrow 0} \partial_k a_{ij}(p) = 0$ is true for any $i,j,k = 1,2$, implying that $a_{ij} \in C^1(\R^2)$. This completes the proof of iii).\hfill $\square$.
\section{Proof of Theorem 1.4}
We start with a general inequality of Caccioppoli type established in [14] dealing with solutions of
\begin{equation}
\divop(\nabla f(\nabla u)) = 0 \quad \text{on } \R^n, n \geqslant 2.
\end{equation}
\begin{lemma} \label{lem3.1}
Let $f \in C^2(\R^n)$ denote a convex function, which means that $D^2 f(p)$ is just positive semidefinite for any $p \in \R^n$. We define
\begin{equation}
\begin{aligned}
\lambda \colon [1, \infty) \to [0, \infty), \quad \lambda(t) &:=
\begin{cases}
t^{-1/2}, & t \ge 2 \\
\frac{1}{\sqrt{2}} (t-1), & 1 \le t \le 2,
\end{cases} \\[1ex]
\Lambda \colon [1, \infty) \to [0, \infty), \quad \Lambda(t) &:=
\begin{cases}
t^{-1/2} (t-1) \ln^2(t), & t > 2 \\
\frac{1}{\sqrt{2}} \ln^2(t), & 1 \le t \le 2.
\end{cases}
\end{aligned}
\tag{3.2}
\end{equation}
Let $\Omega \subset \R^n$ denote a domain and consider a solution $u \in C^2(\Omega)$ of equation (3.1) now on $\Omega$. Then it holds for any $\eta \in C_0^1(\Omega)$ and for all indices $i \in {1,\dots,n}$ ($\Gamma_i := 1 + |\partial_i u|^2$)
\begin{eqnarray*}
\lefteqn{\int_{\Omega} D^2 f(\nabla u)(\nabla \partial_i u, \nabla \partial_i u)\eta^2\lambda(\Gamma_i)dx}\\
& \leq &  c \left( \int_{\text{spt}(\nabla \eta)} D^2 f(\nabla u)(\nabla \partial_i u, \nabla \partial_i u)\eta^2\lambda(\Gamma_i)dx \right)^{1/2}\\
& &\left(\int_{\Omega} D^2 f(\nabla u)(\nabla \eta, \nabla \eta)\Lambda(\Gamma_i) dx\right)^{1/2},
\end{eqnarray*}
where $c$ denotes a constant independent of $\Omega, \eta, u$, and in conclusion\\
\begin{equation}
\int_{\Omega} D^2 f(\nabla u)(\nabla \partial_i u, \nabla \partial_i u)\eta^2\lambda(\Gamma_i)dx \leqslant c \int_{\Omega} D^2 f(\nabla u)(\nabla \eta, \nabla \eta)
\Lambda(\Gamma_i)dx. \tag{3.3}\end{equation}
Moreover, on the right-hand side of (3.3) we may replace $\Lambda(\Gamma_i)$ by $\Lambda(\Gamma)$, where $\Gamma := 1 + |\nabla u|^2$.\end{lemma}
Let us now establish Theorem 1.4: in equation (3.1) we choose the function $f(p) := g(|p|)$, $p \in \R^2$, with $g$ satisfying (1.3), (1.4) and (1.9) and consider a solution $u \in C^2(\R^2)$ with (1.10). For $R>0$ we let $\Omega := B_{2R} := \{x \in \R^2 : |x| < 2R\}$ and define $\eta(x) := \phi(|x|)$, where $0 \leqslant \phi \leqslant 1$, $|\phi'| \leqslant c/R$ on $[0,2R]$ together with \mbox{$\phi = 1$} on $[0,R], \phi = 0$ on $[2R,\infty)$. Fix $i \in {1,2}$  and look at the quantity\\\\
\(A:=\int _{B_{2R}}D^{2}f(\nabla u)(\nabla \eta ,\nabla \eta )\,\Lambda (\Gamma _{i})\,dx\)\\\\
occurring on the right-hand side of (3.3). We have on account of (2.5)
\begin{equation}
D^2 f(p)(q,q) = \frac{g'(|p|)}{|p|}\left[ |q|^2 - \frac{(p \cdot q)^2}{|p|^2} \right] + g''(|p|)\frac{(p \cdot q)^2}{|p|^2}\tag{3.4}
\end{equation}
for points $p, q \in \R^2$, hence we deduce from (3.2) combined with the boundedness of $g'$ and $g''$ (compare (1.4) and (1.9)) that on the set $B_{2R} \cap [|\partial_i u| \leqslant 1]$ the integrand of the quantity $A$ is bounded by $c R^{-2}$, where the positive constant  $c$ is independent of the radius $R$.
\noindent On the set $B_{2R} \cap [|\partial_i u| > 1]$ it holds, again with suitable positive constants $c$ and by using (3.4) together with (1.4), (1.9), (1.10) and (3.2)
\begin{eqnarray*}
\lefteqn{\text{integrand of the term } A} \\
&\leq & D^2 f(\nabla u)(\nabla \eta, \nabla \eta) \ln^2(\Gamma_i)\Gamma_i^{1/2} \\
&\leq & c \Big[ \ln^2(\Gamma_i)g'(|\nabla u|) (|\nabla \eta|^2 - (\nabla u \cdot \nabla \eta)^2 |\nabla u|^{-2}) \Big] \\
&& + \ln^2(\Gamma_i)\Gamma_i^{1/2}g''(|\nabla u|)(\nabla u \cdot \nabla \eta)^2 |\nabla u|^{-2} \Big] \mathrel{:=} c [(\text{I}) + (\text{II})],
\end{eqnarray*}
\noindent
\((\text{II})\leqslant c|\nabla \eta |^{2}\leqslant cR^{-2}\) according to (1.9),\\\\
\((\text{I})\leqslant c\phi ^{\prime }(|x|)^{2}\ln ^{2}(\Gamma )\,|\partial _{t}u|^{2}|\nabla u|^{-2}\leqslant cR^{-2}\) on account of (1.4) and (1.10).\\\\
With (3.3) we altogether obtain
\begin{equation}
\int_{\R^2} D^2 f(\nabla u)(\nabla \partial_i u, \nabla \partial_i u)\lambda(\Gamma_i)dx < \infty,\tag{3.5}
\end{equation}
and from the inequality in front of (3.3) and from (3.5) combined with the previous discussion we deduce
\begin{equation}
\int_{\R^2} D^2 f(\nabla u)(\nabla \partial_i u, \nabla \partial_i u)\lambda(\Gamma_i)dx = 0.\tag{3.6}
\end{equation}
Observing that $D^2 f(p)(q,q) > 0$ for $p, q \in \R^2 - \{0\}$, equation (3.6) implies $\nabla \partial_i u = 0$ für $i=1,2$, which proves our claim. \hfill $\square$
\newpage
% --- Bibliography at the very end ---

\vfill

\noindent Funding.      Open Access funding enabled and organized by Projekt DEAL.\\\\
Data availability.       No additional data are involved.\\\\
Declarations\\
Conflict of interest. On behalf of all authors, the corresponding author states that there is no conflict of interest.

\begin{minipage}[t]{0.48\textwidth}

\textbf{Michael Bildhauer}\\
Saarland University\\
Department of Mathematics\\
P.O. Box 15 11 50\\
66041 Saarbrücken, Germany\\
\texttt{bibi@math.uni-sb.de}

\end{minipage}
\hfill
\begin{minipage}[t]{0.48\textwidth}

\textbf{Martin Fuchs}\\
Saarland University\\
Department of Mathematics\\
P.O. Box 15 11 50\\
66041 Saarbrücken, Germany\\
\texttt{fuchs@math.uni-sb.de}

\end{minipage}
\end{document}